\input amstex
\magnification =1200
\define\spc{strictly pseudoconvex}
\define\sps{strictly plurisubharmonic}
\define\gt{Grauert tube}
\documentstyle {amsppt}
\baselineskip 1.4\baselineskip
\topmatter
\title  Complete hyperbolic Stein manifolds with prescribed automorphism groups\endtitle
\author Su-Jen Kan 
\endauthor 
\address 
Institute of Mathematics, Academia Sinica, Taipei, Taiwan
\endaddress
\email 
{kan\@math.sinica.edu.tw}
\endemail

\NoRunningHeads
\NoRunningHeads
\thanks 2000 {\it Mathematics Subject Classification. } 32C09, 32Q28, 32Q45, 53C24, 58D19. This
research is partially supported by NSC 93-2115-M-001-006. 
\endthanks
\abstract

It is well-known that the automorphism group of a hyperbolic manifold is a Lie group.
Conversely, it is interesting to see whether or not any Lie group could be prescribed as
the automorphism group of certain complex manifold. 

When
the Lie group $G$ is  compact and connected, this problem has been completely
solved by Bedford-Dadok  and  independently by Saerens-Zame  on 1987.  They
have constructed \spc\ bounded domains $\Omega$ such that $Aut(\Omega)=G$. 
For Bedford-Dadok's $\Omega, \ 0\le dim_{\Bbb C}\Omega- dim_{\Bbb R}G\le 1$; 
for generic Saerens-Zame's
$\Omega,dim_{\Bbb C}\Omega \gg dim_{\Bbb R}G$.

J. Winkelmann has  answered  affirmatively to   noncompact connected Lie
groups in recent years. He showed
 there exist Stein complete hyperbolic manifolds $\Omega$ such that $Aut(\Omega)=G$.
In his construction, it is typical that $dim_{\Bbb C}\Omega\gg dim_{\Bbb R}G$.

In this article, we  tackle this problem from a different aspect. We prove
that for any connected Lie group $G$ (compact or noncompact), there exist complete
hyperbolic Stein manifolds $\Omega$ such that $Aut(\Omega)=G$ with $dim_{\Bbb
C}\Omega=dim_{\Bbb R}G.$ Working on a natural complexification of the real-analytic
manifold $G$, our construction of $\Omega$ is geometrically concrete and
elementary in nature.

\endabstract
\endtopmatter

\subheading {1. Introduction}

It is well-known that the automorphism group of a hyperbolic manifold has the structure
of a Lie group. It is natural to ask whether  every Lie group could appear as the
 automorphism group of certain  hyperbolic manifolds or not. This problem sometimes is
called the realization problem: realizing a Lie group as the  automorphism group of
certain complex manifold.

This realization problem has been completely solved 
 by Bedford-Dadok [B-D] and independently by Saerens-Zame [S-Z]  when the
given Lie group is compact and connected. In recent years, generalizing ideas of 
Saerens-Zame, Winkelmann [W] has settled the realizing problem for  any
connected Lie group. However, domains $\Omega$ constructed from this
Saerens-Zame-Winkelmann's approach are typically with $dim_{\Bbb C}\Omega\gg
dim_{\Bbb R}G$.

The strategy of Saerens-Zame-Winkelmann was first to find a domain $D$ on
which $G$ acts by automorphisms, and then perturbed it to a $G$-invariant \spc\ subdomain
in such a way that the additional automorphisms were ruled out by assigning 
CR-invariants to each $G$-orbit on the boundary. To find such a domain $D$ to start
with, Saerens-Zame first embedded the compact Lie group $G$ into the unitary group $U(N_1)$
and then constructed a domain $D$ in $GL(N_1,\Bbb C)\times \Bbb C^{N_2}$ on which $G$ acts
by automorphisms where $N_1$ and $N_2$ are large in general. 
 Due to the above embedding
process, the generic resultant complex manifold
$\Omega$  has  large complex dimension, $dim_{\Bbb
C}\Omega\gg dim_{\Bbb R}G$. Having observed  every Lie algebra is linear and
hence the universal covering of a Lie group could be viewed as linear, Winkelmann
has been able to embed
$\tilde G\hookrightarrow Sp(N_3,\Bbb R)$ and then find a suitable domain $D\subset \Bbb
C^{N_4}$ to start with. Again, this embedding process and the construction of the domain
$D$ has enormously increased the dimension.

The complexification $G_{\Bbb C}$ of a compact Lie group $G$ is Stein with $dim_{\Bbb
C}G_{\Bbb C}=dim_{\Bbb R}G.$
Starting from  domains in $G_{\Bbb C}$, Bedford-Dadok were able to give a more concrete
construction. They found bounded \spc\ domains
$\Omega\subset G_{\Bbb C}$ or $\Omega\subset G_{\Bbb C}\times\Bbb C$ such that
$Aut(\Omega)=G$.

A natural attempt is to generalize Bedford-Dadok's approach to noncompact Lie groups.
Unfortunately,  it is not easy to give a canonical complexification of a generic
noncompact  Lie group such that $dim_{\Bbb
C}G_{\Bbb C}=dim_{\Bbb R}G.$  In this article  we
consider a special kind of complexifications of real-analytic manifolds to
resolve this difficulty. 

A Lie group $G$ equipped with a
left-invariant metric $g$ is naturally a homogeneous space with $Isom (G)\simeq
L(G)\cdot K$ where
$K$ is the isotropy group and $L(G)$ are left multiplications.
 We provide the
tangent bundle of the Riemannian manifold $(G,g)$ with
 a complex structure   in such a way that 
all the leaves of the  Riemannian foliation are holomorphic curves.
The  disk bundle of radius $r$ is called  a \gt \ $T^rG$. 

Generalizing  the rigidity result, $Aut(T^rG)=Isom (G)$, proved in [K]  we are able to
dominate automorphism groups of certain  domains $D\subset T^rG$.  We prove it at Theorem
3.5 and call it the  subrigidity $Aut(D)<Isom(G)$.

Starting from such a domain $D$ our strategy is to destroy, through certain $G$-invariant
perturbations, the additional automorphisms coming from the isotropy group.
Though the objects we deal with are not even relatively compact, the holomorphic 
extension property needed here is not hard to handle due to the 
special complex structure adopted here.
In fact, most mappings we deal with are bundle mappings which automatically
extend over the boundary. 
The main result of this article is

\proclaim{Theorem }
Let $G$ be a  connected Lie group of dimension $n\ge 2$. Then there exist   complete
hyperbolic Stein manifolds
 $\Omega, dim_{\Bbb C}\Omega=n$,  such that $Aut(\Omega)=G$.
\endproclaim 
The dimensional
condition
$n\ge 2$ has to be added here since the main idea we'll use  follows from the
rigidity arguments  of \gt s while \gt s and their perturbations are biholomorphic to the
unit disc, by the Riemann mapping theorem, when the center  is of dimension one. 
In the rigidity arguments of  \gt s, it is also necessary to assume the Riemannian
manifold is connected. At this moment, I could not see a way to release  the
connectedness of the Lie groups.
We remark while Saerens-Zame's result works for compact disconnected  Lie group as well,
the resultant complex manifold $\Omega$ is typically with $dim_{\Bbb
C}\Omega\gg dim_{\Bbb R}G$.

We emphasize that  methods developed in this article work simultaneously for 
compact and noncompact Lie groups; the complex dimension of the resultant  complex
manifold 
$\Omega$ is equal to the real dimension of $G$; and
  $\Omega$ thus obtained are
constructed in a geometrically concrete and  nature way. 

The organization of this article is the following. In $\S2$ we briefly review 
terminologies concerning \gt s and prove the existence of Stein \gt s. Generalizing 
notations about the rigidity of \gt s, we prove subrigidity characterization of
certain domains in 
$\S3$. In $\S4$, specific domains and perturbations are constructed explicitly.
We perturb domains in the tangent bundle of  a connected Lie group in an invariant
way such that extra  symmetry on each fiber would be eliminated. By constructed such kind
of domains in a fairly explicitly way, the realization of a connected Lie group as an
automorphism group   follows from the subrigidity derived on
$\S3$.

This realization problem was first noticed by the author and also mentioned to the
author by  Peter Heinzner  in September 2002 when I was visiting
Ruhr-Universit\"at, Bochum Germany. Over years I was trapped in an effort to show every
\gt\ is Stein. Among other things, I am very grateful to 
Peter Heinzner for bringing my attention to  [HHK] and [D-G],
which have helped settling down the Steiness part  needed in this article, during my
recent visit to Bochum in August 2004. I would also like to thank J\"org Winkelmann for
many helpful discussions about this project during the Hayama Symposium  2002.

\

\subheading{2. Existence of Stein Grauert tubes}

 Throughout this article,   $(M,g)$ is assumed to be  a connected
real-analytic Riemannian manifold of dimension $\ge 2$.  There exists a complex
structure, the {\it adapted complex structure}, in a domain $\Omega(M) \subset TM$ which
turns leaves of the Riemannian foliation on
$\Omega(M)\backslash M$ into holomorphic curves. Let's denote by $\Omega(M)$  the
maximal domain in which the adapted complex structure exists. With respect to this 
complex structure, the length square function 
$\rho (x,v):=|v|^2$ is
\sps\  and satisfies the complex homogeneous Monge-Amp\`ere equation $(d d^c\rho)^n=0$
on
$\Omega (M)\backslash M$ with the initial condition   $\rho_{i\bar
j}|_{\sssize M}=\frac 12 g_{i j}.$ 
The {\it
\gt} 
$$T^rM=\{(x,v):x\in M, v\in T_xM, |v|<r\}=\{\rho^{-1}[0,r^2)\}$$ is the collection of
tangent vectors of length less than $r$   equipped with the  adapted complex 
structure; $M$ is the  center and $r$ is the radius. In general the maximal
domain $\Omega(M)$ is not a \gt, $\Omega(M)$ is a \gt\ if and only if $M$ is a symmetric
space of rank one.
 There is a natural antiholomorphic involution $\sigma$ fixing
every point of  $M$,
$$\sigma: \Omega (M)\to  \Omega (M),\;\; (x,v)\to  (x,-v).\tag 2.1$$

  For two isometric real-analytic Riemannian manifolds $(M,g)$ and $(N,\kappa)$, the
nature of the adapted complex structure will assure the biholomorphic equivalence of
$\Omega(M)$ and
$\Omega(N)$ and the biholomorphic equivalence of $T^rM$ and $T^rN$ as well. Given
an isometry $h$ of $(M,g)$, the differential $dh$ acts as a biholomorphism on
$\Omega(M)$;
$$dh:\Omega(M)\to \Omega(M),\;\;
dh(x,v)= (h(x), h_*v).\tag 2.2$$ 
The notation $F_p^r$ is reserved for the fiber passing $p\in M$,
$$F^r_p:=\{(p,v):v\in T_pM, |v|<r\}.\tag 2.3$$

Since the complex structure we consider here is a local object, for each real-analytic
Riemannian manifold, there exists an 
$r_{max}(M)\ge 0$, the maximal radius such that the adapted complex structure could be
 well-defined   on
$T^{r_{max}}M$. For compact, co-compact or homogeneous  $M$,  $r_{max}(M)>0.$  For any
$r<r_{max}(M)$, the
complex manifold $T^rM\subset T^{r_{max}}M$ with \spc\ boundary
defined by $\{\rho=r^2\}$. 

It is shown in Theorem 5.2 of [K] that any \gt\ over homogeneous space of radius
$r<r_{max}$ is complete hyperbolic. Another
natural question  to ask about is whether or not a \gt\ is  Stein. 

 Recall from [D-G] that an {\it unramified Riemannian domain} $D$ over a Stein manifold
$X$  is a complex manifold $D$ together with a locally homeomorphic holomorphic map from
$D$ to $X$ and that a pseudoconvex unramified Riemannian domain over a Stein manifold
is itself Stein. It is clear that any pseudoconvex domain in a Stein manifold is 
Stein since we may take  the imbedding as the  holomorphic map.

\proclaim{Proposition 2.1}
$T^rM$ is a Stein manifold for all $r\le r_{max}(M)$ if $M$ 
  is one of the following:  compact; co-compact;  homogeneous with nonnegative
curvature.
\endproclaim
\demo{Proof}
 For compact  $M$, the \gt\ $T^rM,r\le r_{max}$ is exhausted by the \sps\ function
$-\log (r^2-\rho)$, hence is  Stein.  

 When $M$ is co-compact,  $\hat M=M/\Gamma$ is compact for some discrete
subgroup $\Gamma<Isom(M)$, then the \gt\ $T^r  M$ 
 is simply a covering  of the Stein \gt \ $T^r\hat M$. Thus 
$T^rM, r\le r_{max}(M),$
is Stein.

 Following the decomposition theorem( Theorem 7.1*, 
[C-G]) of Cheeger-Gromoll  a homogeneous manifold $M$ of nonnegative curvature may be
written as the product  $\Bbb R^k\times  M^*$ where $ M^*$ is a compact
homogeneous space of nonnegative curvature. Since the metric on $M$ is the product metric
from 
$\Bbb R^k$ and $ M^*$,
it is clear that $r_{max}(M)=r_{max}( M^*)$ and
$T^r(M)\subset \Bbb C^k\times T^{r_{max}}( M^*), \,\ \forall r\le r_{max}.$
For any $r\le r_{max}$, $T^r(M)$ has pseudoconvex boundary since the boundary is 
locally exhausted by the \sps \ function $-\log (r^2-\rho)$.
Being a  pseudoconvex  domain in the Stein manifold $\Bbb C^k\times
T^{r_{max}}( M^*)$,  $T^{r_{max}}M$ is Stein. \qed
\enddemo
Though there is no example disproving the Steiness of any \gt, the above three
categories are the only complete classes that we are sure about the Steiness of any \gt s
over. For some very special  kind  of homogeneous spaces $M=G/K$, Halverscheid-Iannuzzi
[H-I] were able to construct a polar map from  $TM$ to the Stein manifold $G_{\Bbb
C}/i(K_{\Bbb C}).$ They show the maximal domain over a Heisenberg group is neither
holomorphically separable nor holomorphically convex. However, there are some
facts hiding behind: Grauert tubes constructed over their spaces
 are  Stein. Their polar map is an unramified one therefore, \gt s are
pseudoconvex Riemann domains over the Stein manifold  $G_{\Bbb
C}/i(K_{\Bbb C}),$ they are Stein (c.f. [D-G]).

 A
complex manifold $X$ is  a  {\it complexification} of $M$ if  $M\subset X$ as a maximal
totally real submanifold  $dim_{\Bbb C}X=dim_{\Bbb R}M$. Utilizing his solution to
the Levi problem, Grauert has proved  there always exists a Stein complexification of a
real-analytic manifold.

Though any two complexifications $(X_1,J_1)$ and $(X_2,J_2)$ of a real-analytic manifold
$M$ are locally biholomorphic near $M$, i.e., there exist a neighborhood $U_1$ of $M$
in $X_1$ and a neighborhood $U_2$ of $M$ in $X_2$ such that $U_1$ and $U_2$ are
biholomorphically equivalent, we still can't conclude there are some Stein \gt s of
small radii since it does not seem clear how to control the radius uniformly when the
manifold is  not relatively compact.

In [HHK], the authors have put on an extra piece of information, say there is some kind of
$G$-action on $M$. They ask for the question whether the
 Stein manifold in Grauert's complexification could still be chosen to be
$G$-invariant or not. Analyzing the real-analytic slices of certain categorical
quotients, they  conclude the following:
 \proclaim{Theorem(H-H-K)} Let $X$ be a 
$G$-complexification of a real-analytic manifold $M$ where $G$ is a connected Lie group
 acting properly on  $M$ as real-analytic diffeomorphisms. 
Then there exists a 
$G$-invariant Stein neighborhood of $M$ in $X$.
\endproclaim 

A group action $G$ on $M$ is  {\it proper} if  the inverse image of a
compact set is compact. A  complexification $X$  of $M$ is called a {\it
$G$-complexification}  if the
$G$-action on $M$ extends  to a holomorphic $G$-action on $X$.

Though it is not clear whether  a
generic \gt\ is Stein or not, a direct application of the above theorem would
assert the existence of Stein \gt s of small radii provided the centers have possessed
some transitivity property. This would be enough for the purpose of solving the
realization problem.

\proclaim{Proposition 2.2}
Let $(M,g)$ be a homogeneous space. Then there exist $\epsilon >0$ such that the \gt \
$T^{\epsilon}M$ is Stein and complete hyperbolic.
\endproclaim
\demo{Proof}
Let $G=Isom_0 (M)$, then $G$ acts properly on $M$ and acts
holomorphically on the complexification  $T^{r_{max}}M$ of $M$. By [HHK], there 
exists a $G$-invariant Stein neighborhood $U$ of $M$ in $T^{r_{max}}M$. For $p\in M$ we
may assume the fiber $F^{\epsilon}_p$,  tangent vectors at $p$ of length less
than
$\epsilon$, is contained in  $U$. By the transitivity of the $G$-action, $G\cdot
F^{\epsilon}_p$ is the \gt\ $T^{\epsilon}M.$ The Steiness follows since 
$T^{\epsilon}M$ is now a \spc\ domain in the Stein manifold $U$.\qed
\enddemo

\

\

\subheading{3. Subrigidity of  domains}

Let $Aut$ denote the automorphism group and $Isom$ denote the isometry group; $Aut_0$
and $Isom_0$ denote the corresponding identity components. For  any  real-analytic
Riemannian manifold
$(M,g)$ such that 
$T^rM$ is not covered by the ball,
 the following  rigidity results were proved  in  [K].
 \roster
\item $Aut_0(T^rM)=Isom_0(M)$;
\item furthermore, if $M$ is homogeneous then $Aut(T^rM)=Isom(M)$.
\endroster

We remark here that whenever the rigidity is mentioned, we always assume $\dim M>1$.
When $\dim M=1$, the \gt\ $T^rM$ is the disc in $\Bbb C$. 
 One essential feature of \gt s is the symmetry on each fiber, each fiber is
a disk bundle. In this section, we'll consider domains in a more
general setting, perfect symmetry on fibers won't be asked for any more. One crucial step in
proving the rigidity  is the Theorem 4.1 in [K] which has characterized the isometry
group of $M$. We will show  this kind of characterization still works and   
$Aut(D)$ is  dominated by
$Isom(M)$ provided there is some transitivity on $M$.
Arguments will follow  the  spirits and methods developed in [K]. 
However, 
an extra piece of assumptions on the target \gt\ is necessary; we assume  \gt s are
Stein while the existence has been guaranteed in the last section. 
\proclaim{Lemma 3.1}
Let $T^rM$ be a Stein \gt, $D\subset T^rM$ be a connected domain containing  $M$. If
$f\in Aut (D) , f(M)=M$, then $f=du$ for some isometry $u$ from $(M,g)$ to
$(M,\kappa)$ and
$f$ could be extended
 over the boundary $\partial D$.

\endproclaim

\demo{Proof}
From the construction of \gt s, it is clear that $M$ is a maximal totally real
submanifold of $D$. Let $\rho$ denote the length square function on $T^rM$, then 
$\rho\cdot f^{-1}$ is a
\sps \ function on $D$. The Riemannian metric $g$ is induced from the K\"ahler form 
${i}\partial\bar\partial\rho$, $2\rho_{i\bar j}|_{\sssize
M}= g_{i j}.$  Similarly, the  K\"ahler form 
${i}\partial\bar\partial(\rho\cdot f^{-1})$ has induced a Riemannian metric $k$ on
$M$. Let $u=f|_{M}:(M,g)\to (M,\kappa)$; $u$  is an isometry. 
By the nature of the adapted complex structure, the differential $du$ is a
biholomorphic mapping from the \gt\ $T^r(M,g)$ to the \gt\ $T^r(M,\kappa).$ 

Let $D_1=D\cap T^r(M,\kappa)\hookrightarrow \Bbb C^N$ for some $N$. This could be
achieved since $T^r(M,g)$ is Stein and $T^r(M,\kappa)$ is biholomorphic to $T^r(M,g)$. Let
$D_2=(du)^{-1}(D_1)\subset T^r(M,g)$. Let
$D_3=D\cap D_2$ then $D_3$ is a
domain in
$T^rM$ which has contained
$M$ as a maximal totally real submanifold. Restricted to $D_3$, 
  $f$ and $du$ could be viewed  as   holomorphic functions from $D_3$ to $\Bbb C^N$.
Since $f=du$ at the maximal totally real submanifold $M$, by the identity principle,
$f\equiv du$ on the subdomain $D_3$ and hence  $f\equiv du$ on $D$. For $(x,v)\in D,
f(x,v)=du(x,v)=(u(x),u_*v))=(f(x),u_*v)$ is a bundle map acting fiberwise. It is clear
that 
$f$ could be extended  over the boundary of $D$.
\qed
\enddemo
 
Recall $\Omega(M)$ is the maximal domain in $TM$ such that the adapted complex structure
is defined; $\sigma$ is the natural antiholomorphic involution in $\Omega(M)$ and  
the norm-square function $\rho$ is 
\sps\ in $\Omega(M)$. Sibony (Theorem 3, [S]) has asserted the hyperbolicity of any 
complex manifold equipped with a bounded \sps\ function. For domains in $\Omega(M)$ to
be hyperbolic, the only thing we got to take care is the vertical direction. 

One essential
feature of Grauert tubes is that each fiber, which is diffeomorphic to a real
ball, has perfect symmetry. We have tried to release this symmetry on the fiber 
to get some kind of rigidity result. However, it seems to us  certain kind of
transitivity on  the domain is necessary for our purpose. We make the following
definition.
\proclaim{Definition 3.2}  A  domain  $D\subset\Omega(M)$ is  
$G$-homogeneous if:
\roster 
\item there exists a connected subgroup
$G\subset Isom(M)$  acting transitively on
$M$;
\item there exists a bounded open subset $ F_p\subset T_pM$ such that 
$D=G\cdot F_p$.\endroster 
\endproclaim

With this extra piece of symmetry,  Lemma 3.1 has the following refined form:
\proclaim{Proposition  3.3} 

Let $D$ be a  $G$-homogeneous domain in $\Omega(M)$.
If $f \in Aut (D),f(M)=M$, then
$f=du$ for some $u\in Isom(M)$.
\endproclaim
\demo{Proof}
Since $M$ is homogeneous, there exists small $\epsilon$ such that $T^{\epsilon}M$ is
Stein. The restriction map of $f$ to $D\cap T^{\epsilon}M$ has possessed the bundle map
property developed in Lemma 3.1. Thus $f=du$, for some isometry $u$ from the Riemannian
manifold $(M,g)$ to the Riemannian manifold $(M,\kappa)$. By the assumption on $D$ 
there exist $G\subset Isom(M)$ and $p\in M$ such that 
$D=G\cdot F_p$. 

Let $q=u(p)=h(p)$ for some $h\in G$ and  let $\{e_1,\cdots,e_n\}$
be an orthonormal basis of $F_p$ with respect to the metric $g$. Since $h$ is an
isometry,
$\{(e_1,\cdots,e_n)A\}$ is an orthonormal basis of $F_q$ with respect to the metric $g$
where $A\in O(n)$. On the other hand, $\{u_*e_1,\cdots,u_*e_n\}$ form an orthonormal
basis of $F_q$ with respect to the metric $\kappa,$ thus, there exists a $B\in GL(n,\Bbb
R)$ such that $$(u_*e_1,\cdots,u_*e_n)=(e_1,\cdots,e_n)AB.$$ Let's denote the matrix
$AB=C=(C_1,\cdots, C_n)$ and $e=(e_1,\cdots,e_n)$ then 
$$\delta_{ij}=<u_*e_i,u_*e_j>_{\kappa}=<e C_i,eC_j>_g=C_i\cdot C_j.$$
This shows  $C$ is a matrix with orthonormal columns, thus $C\in O(n)$ and $B\in
O(n)$ as well.

The orthonormal basis $\{u_*e_1,\cdots,u_*e_n\}$ of the metric $\kappa$ come from an
orthogonal transformation of the orthonormal basis $\{e_1,\cdots,e_n\}$ of the metric
$g$. Thus, $\kappa=g$ and $u\in Isom(M,g)$.
\qed
\enddemo

A $G$-homogeneous domain  is hyperbolic and hence its automorphism group is
a Lie group. Let $D$
be a $\sigma$-invariant \spc\ $G$-homogeneous domain in $\Omega(M)$. 
Following notations used  in $\S6$ of [K], we set $\widehat {Aut}(D):=
Aut(D)\cup\sigma\cdot Aut(D)$ and denote
$\hat\Cal L$ as the Lie algebra of $\widehat {Aut}(D)$. Modifying the proofs of
Lemma 6.1, Lemma 6.2 and Lemma 6.3 in [K], it could be shown that for every $\xi\in\Cal
L$, the Lie algebra of $Aut(D)$, the vector field associated to $\xi$ is tangent to $M$.
We conclude

\proclaim{Proposition 3.4}
 Let $D$ be a $\sigma$-invariant \spc\  $G$-homogeneous domain in $\Omega(M)$. If
$D$    is not biholomorphic to the ball. Then $Aut_0(D)\subset Isom_0(M).$
\endproclaim

Furthermore, following [K] again, we are going to
prove the subrigidity $Aut(D)\subset Isom(M)$. For the rest of this
section, we assume $D$ is a \spc\ $\sigma$-invariant  $G$-homogeneous  domain in
$\Omega(M)$, $D$ is not biholomorphic to the ball. Proofs go exactly the same way as in
$\S7$ of [K]; a brief explanation would be given in the following. 

(I). First of all, a \spc\  $G$-homogeneous  domain in $\Omega(M)$ is complete
hyperbolic since the boundary behavior is dominated by the \spc ity and the horizontal
direction is determined by the transitivity of the $G$-action. 

(II). Secondly, the following is clear
 from the $G$-homogeneity and Prop. 3.4.
$$G\subset Aut_0(D)\subset Isom_0(M).$$
Given $f\in Aut (D)$, $D$ is also a  $f\cdot G\cdot f^{-1}$--homogeneous domain
centered at $N=f(M)$. By the homogeneity of the $G$-action,
$$\aligned & N=f(M)=f(G\cdot p)=f(Aut_0(D)\cdot p)=Aut_0(D)\cdot f(p);\\
& \pi (N)=\pi (Aut_0(D)\cdot f(p))=Aut_0(D)\cdot \pi (f(p))=M.
\endaligned$$

Thus, the argument in Lemma 7.1 of [K] goes through and $f(M)\cap T_pM\ne\emptyset$ for
any $p\in M$. 

(III). Since $D$ is  complete hyperbolic and $G$-homogeneous,  arguments in Prop.7.2
of [K] could be transplanted here.
Thus, the index of $G$ in
$Aut(D)$ is finite.

(IV). Let $\tau=f\cdot\sigma\cdot f^{-1}$ be the associated antiholomorphic
involution with respect to the fixed point set $N=f(M)$, then the identity relation  at
(7.2) and the Proposition 7.3 of [K] work as well. Thus, there exists an odd integer $k$
such that
$(\sigma\cdot\tau)^k=id.$

We conclude this section by the following theorem.
\proclaim{Theorem 3.5}

Let $D$ be a $\sigma$-invariant \spc\  $G$-homogeneous  domain in $\Omega(M).$
Then either  $Aut(D)\subset Isom(M)$ or $D$ is the ball.
\endproclaim
\demo{Proof}
As stated in (III),  $Aut(D)/G=\{g_jG: g_j\in Aut(D), j=1,\dots, k\}$.
Then $\psi(z)=\sum_{j=1}^k \rho(g_j(z))$ is an
$\widehat {Aut}(D)$-invariant strictly plurisubharmonic
nonnegative function in $D$ where $\rho$ is the length square function in $D$. As $G$
acts transitively on
$M$, the tangent space $T_z(D)$ could be decomposed as, for any $z\in D$, 
$$T_z(D)=T_z(G \cdot z)+T_z(T_{\pi(z)}M\cap D)$$ where $\pi$ is the fiber projection
$\pi (x,v)=x,\; \forall v\in T_xM$.

 Since $\psi$ is constant in $G
\cdot z$, every critical point of the function $f_z:= \psi|_{T_{\pi(z)}M\cap D}$ is a
critical point of $\psi$ and every critical point of $\psi$ occurs at the critical
points of the functions $f_z$.

 As $\psi$ is strictly plurisubharmonic, the above
decomposition implies that the real Hessian of $f_z$ is positive definite on the tangent
space 
$T_z(T_{\pi(z)}M\cap D).$ Since $f_z$ is proper on the  fiber, it follows
that there is exactly one critical point of $f_z$ which turns out to the  minimal point.
Since $\psi\cdot
\sigma=\psi$, the minimum of $f_z$ occurs at $\pi (z).$ That is to say that the set of critical
points of $\psi$ is exactly  $M$.

Let $f\in Aut(D), f(M)=N$, then $N$ is the critical point set of $\psi$ since $\psi$ is
$Aut(D)$-invariant. We conclude that $N=M$ and  $f\in Isom(M)$ following from Prop.3.3. 
The inclusion $Aut(D)\subset Isom(M)$ is concluded.\qed
\enddemo

\

\subheading{4. Realizing a connected Lie group as an automorphism group}

Let $G$ be a connected real Lie group with Lie algebra $\Cal G$. Given a positive
definite inner product $<\;,\;>$ on $T_eG=\Cal G$, we may endow  $G$ with the associated 
 left invariant Riemannian metric $g$. Every Lie
group is real-analytic, since it is
locally diffeomorphic to the Lie algebra $\Cal G$ through the exponential map. Thus
$(G,g)$ is a real-analytic Riemannian manifold. Furthermore, $(G,g)$ is a homogeneous
space with trivial tangent bundle $G\times T_eG$.

It is clear
that the left translation group $L(G)$ is a subgroup of $Isom(G,g)$ and $L(G)$ acts
transitively on $G$. Furthermore, the transitivity implies the  diffeomorphic
equivalence of $G$ with
$Isom(G,g)/K$ and with $Isom_0(G,g)/K_0$ where 
$$K=\{h\in Isom(G,g):h\cdot e=e\}\tag 4.1$$ is the
isotropy group at $e\in G$ and $K_0=K\cap Isom_0(G,g)$. The
following two equations are  immediate.
$$Isom(G,g)\simeq L(G)\cdot K,\;\; Isom_0(G,g)\simeq L(G)\cdot K_0.\tag 4.2$$

From now on, we denote the homogeneous Riemannian manifold $(G,g)$ as $(M,g)$ where $G$
is a connected Lie group of dimension $n\ge 2$ and $g$ is a left invariant real-analytic
metric on $G$. We would like to construct a 
$G$-homogeneous
\spc\  domain  and then further perturb the domain  to
eliminate additional automorphisms such that the automorphism group of the resultant
domain
 is  $G$. We emphasize  that our method work for both
compact and noncompact Lie groups as long as the Lie group is connected. The dimensional
condition
$n\ge 2$ has to be added here since the main idea we used here follows from the
rigidity argument  of \gt s while the rigidity of \gt s fails when the center is of dimension
one. 

Let $(M,g)=(G,g)$ be as above. Since it is homogeneous, there exist $\delta>0$ such
that the \gt\ $T^{\delta}M$ is  Stein  and the \gt\ $T^{2\delta}M$ still exists. Let
$\rho$ denote the length square function; $\sigma$ denote the natural antiholomorphic
map in $T^{2\delta}M$ and 
$F^r_p$ denote tangent vectors at $p\in M$ of length less than $r$,
$F^r_p:=\{(p,v):v\in T_pM, |v|<r\}.$ 

Let $\{e_1,\cdots,e_n\}$ be an
orthonormal basis--with respect to the metric $g$-- of $F^{2\delta}_p$, $K$ be the
isotropy group of $(M,g)$ at $p\in M$. Then $K\subset O(n)$ when we view
$F^{2\delta}_p$ as a subspace of the real vector space generated by the orthonormal basis
$\{e_1,\cdots,e_n\}$. 

Let $a_j=\delta e_j\in  F^{2\delta}_p, j=1,\cdots,n.$ Choose a ball
$B_{\epsilon}\subset F^{2\delta}_p$ centered at
$a_1$ of  radius $\epsilon\ll\delta$ and a family of orthogonal transformations
$f_j\in O(n)$ such that $ f_j(a_1)=a_j,j=2,\cdots,n.$ 

Let $\eta_1(x)=(x-a_1)^{4l},l\ge 1,$ be a real-analytic function on 
$B_{\epsilon}$. Convoluting with some cut-off function, we may assume $\eta_1$ has
compact support
$C_1\subset B_{\epsilon}$. 

Denoting $f_1=id$, we define a function $\eta_2$ on $F_p^{2\delta}$ as follows:
$$\eta_2(f_j(x))=\eta_2(\sigma (f_j(x)))=\frac 1j\eta_1(x),\;\forall x\in F_p^{2\delta},
j=1,\cdots,n.\tag 4.3$$

The ball $B_{\epsilon}$ could be arranged so small that  $\{\sigma^i
(f_j(B_{\epsilon})): i=0,1; j=1,\cdots,n\}$  are pairwise disjoint.
It is clear that $\eta_2$ is a $\sigma$-invariant function in $F_p^{2\delta}$
 with compact support
$C_2$ contained in
$\cup_{j=1}^n (f_j(B_{\epsilon})\cup \sigma (f_j (B_{\epsilon}))).$

Since the tangent bundle of a Lie group is trivial, $TM=G\times T_pM$.
We define a real-analytic functions $ \eta$ in $T^{2\delta}M=G\times F_p^{2\delta}$ by
setting 
$$\eta(h\cdot x)=\eta_2(x),\;\forall h\in G, x\in F_p^{2\delta}.\tag 4.4$$
 
Clearly, $ \eta$ is a nonnegative $G$-invariant and $\sigma$-invariant real-analytic
function on
$T^{2\delta}M$.
Recall $\rho$ is the length square function of the tangent vectors which is
$G$-invariant and 
\sps\  and
$T^{\delta}M=\{\rho^{-1}([0,\delta^2))\}$. 
Define the function, for some $0<\epsilon'\ll 1$, $$\hat \rho=\rho
-\delta^2+\epsilon'\eta\tag 4.5$$  Shrinking $\epsilon'$ if necessary such that  all the 
0-th, first and second order derivatives of $\epsilon'\eta$ are well under control.
Thus, we may assume the  function $\hat \rho$
 is \sps. It is $G$ and $\sigma$
invariant since both $\rho$ and $\eta$ are.

Let $D$ be a domain defined by the function $\hat\rho,$
$$D:=\{z\in T^{2\delta}M;\hat\rho (z)<0\}.\tag 4.6$$ 

It is clear that $D$ is a $\sigma$-invariant \spc\  $G$-homogeneous domain in
$ T^{\delta}M$. In order to apply Theorem 3.5 to conclude the subrigidity, we
need to show the domain $D$ is  not biholomorphic to the ball. Some background on
Chern-Moser normal form is needed here.

In the fundamental paper [C-M], Chern and Moser  have
associated to every \spc\ point $p$ in a hypersurface $H$ a family of local
 invariants, namely  a neighborhood $U_p$ of $p$ in 
$H$ is biholomorphically equivalent to a neighborhood $V_q$ of $q$ in a hypersurface 
$S$ if and only if the associated families of invariants at $p$ and at $q$ are the
same. These invariants are  given by the coefficients of certain  normal form
of the defining function which we briefly explain in the following. Let $\psi$
be a local defining function of the $(2n-1)$-dimensional  hypersurface $H$ near the 
point $p$, say inside a coordinate chart. 
Since $p$ is \spc\ the Levi form is positive definite and the defining function could
be transformed  through some linear translation and  proper
holomorphic transformations to the following:
$$v=|z|^2+ F(z,\bar z, w)\tag 4.7$$
where $z=(z_1,\cdots,z_{n-1})\in \Bbb C^{n-1}, w=u+iv\in \Bbb C$. Using the 
transformations
$z^*=z+f(z,w),\,\, w^*=w+g(z,w)$
Chern and Moser have simplified (4.7) to  
$$v^*=|z^*|^2+ N_{24}(u^*)z^{*2}\bar z^{*4} + N_{42}(u^*)z^{*4}\bar z^{*2}+ \sum_{j,k\ge
2,j+k\ge7}N_{jk}(u^*)z^{*j}\bar z^{*k}\tag 4.8$$
Furthermore this transformation, and hence all of the coefficients $N_{jk}$,  is made
unique when certain normalizations on $f$ and $g$ are made.

For the hyperquadric $Q=\{v=|z|^2\}$, the unbounded model of the sphere, all
invariants $N_{jk}$ vanish.

\proclaim{Lemma 4.1} $D$ is not biholomorphic to the
ball.
\endproclaim
\demo{Proof} It is standard
that given any $n$-dimensional totally real closed submanifold $X\subset B^n$ there
exists $\varphi\in Aut(B^n)$ such that $\varphi (X)= \Bbb R^n\cap B^n$.

 Suppose there exists a biholomorphic map
$f:D\to B^n$, then $f(M)$ is a maximal  totally real closed submanifold of
$B^n$ since $M$ is such kind of submanifold of $D$.
Without loss of generality, we may assume $f(M)=\Bbb R^n\cap  B^n$.
Adopting the argument used in the proof of Lemma 3.1, we see $f$ is a bundle map, $f=du$
for some isometry
$u$ from  $(M,g)$ to $(\Bbb R^n\cap B^n,g^*).$ 

Thus, the biholomorphic map $f$
could be extended holomorphically over the boundary $\partial D$ and thus 
 the Chern-Moser normal form
of $\hat\rho$ at any boundary point of $D$ has the local expression $v=|z|^2$.

 However,
the defining function  near $a_1\in \partial D$ and the defining function near generic
boundary points differ by some  $4l$-order
terms. In the construction of the
normal form, we see the orders won't decrease.
Thus the Chern-Moser normal form at
$a_1$ won't be the same as the Chern-Moser normal form at generic boundary points
which we assume to be $v=|z|^2$. A contradiction.\qed
\enddemo

It is clear $G$ is contained in the automorphism group of $D$ since $D$ is
$G$-invariant.
 Theorem 3.5 along with (4.2) implies $$L(G)<Aut(D)<Isom (M)\simeq L(G)\cdot K.\tag
4.9$$ The last step is to eliminate those automorphisms coming form the isotropy group.
\proclaim{Theorem 4.2}$D$ is a complete hyperbolic Stein manifold with
$Aut(D)=G$.
\endproclaim
\demo{Proof} $D$ is complete hyperbolic since it is   $G$-homogeneous with
\spc\ boundary. It is Stein since $D$ is a \spc\ domain in the Stein manifold
$T^{\delta}M$.

It remains to show if $h\in K\cap Aut(D)$, then $h$ is the identity
map.  By the construction of $\eta_2$ at (4.3) and hence the construction of
the defining function $\hat\rho$, we see that there exist neighborhoods $U_j$ of $a_j$
on the hypersurface $\partial D\cap F_p^{\delta}$ such that points in $U_i-a_i$ all have
norms different from  norms of points in $U_j-a_j$ for $i\ne j$. Besides, every
point in
$U_i-a_i$ has norm $<\delta$ and $|a_i|=\delta,\; i=1,\cdots,n$.

Since $K$ is a subgroup of $Isom(M), h$ is a linear and norm-preserving
isomorphism of $F^{\delta}_p\cap D.$ In fact, $h$ is the restriction of an orthogonal
transformation in $F^{\delta}_p$. Thus, $h$ is a norm preserving map from $\partial
D\cap F_p^{\delta}$ to $\partial D\cap F_p^{\delta}$ and the only possibility for $h$
is either $h(a_j)=a_j$ or $h(a_j)=\sigma (a_j)=-a_j,\; j=1,\cdots,n.$

Recall that  $a_j=\delta e_j$, the linearity of $h$ implies either $h(e_j)=e_j$
or $h(e_j)=-e_j,\forall j$.  Thus, $h$ is either the identity map or the negative
identity, in other words, $h=id$ or $h=\sigma$. The second case is not possible since
then $h$ is antiholomorphic rather than holomorphic. We conclude $h=id, K=id$ and
$Aut(D)=L(G)\simeq G$.\qed

\enddemo

\enddocument